\documentclass[10pt,a4paper]{amsart}
\usepackage{graphicx,latexsym,amssymb}
\usepackage[all]{xy}
\CompileMatrices
\theoremstyle{plain}
\newtheorem{theorem*}{Theorem}
\newtheorem{theorem}{Theorem}[section]
\newtheorem{lemma}[theorem]{Lemma}
\newtheorem{sublemma}{Sublemma}
\newtheorem{proposition}[theorem]{Proposition}
\newtheorem{corollary}[theorem]{Corollary}
\newtheorem{corollary*}[theorem*]{Corollary}
\theoremstyle{definition}
\newtheorem{definition}[theorem]{Definition}
\newtheorem{example}[theorem]{Example}
\theoremstyle{remark}
\newtheorem{remark}[theorem]{Remark}

\numberwithin{equation}{section}
\numberwithin{figure}{section}
\font\nb=msbm10
\font\nbi=msbm8 at 7pt
\def \R{\hbox{\nb R}}
\def \N{\hbox{\nb N}}

\def \Z{\hbox{\nb  Z}}
\def \Zi{\hbox{\nbi Z}}
\def \Q{\hbox{\nb Q}}
\def \Qi{\hbox{\nbi Q}}
\def \C{\hbox{\nb  C}}

\begin{document}
\title[]
{Quadratic functions on torsion groups}
\date{\today}
\author[F. Deloup]{Florian Deloup}
\author[G. Massuyeau]{Gw\'ena\"el Massuyeau}

\address{Florian Deloup,
Einstein Institute of Mathematics,
Edmond J. Safra Campus, Givat Ram,
The Hebrew University of Jerusalem,
91904 Jerusalem, Israel, {\it{and}}}
\address{Laboratoire Emile Picard,
UMR 5580 CNRS/Universit\'e Paul Sabatier,
118 route de Narbonne,
31062 Toulouse, France.
}
\email{deloup@math.huji.ac.il}

\address{Gw\'ena\"el Massuyeau,
Laboratoire Jean Leray,
 UMR 6629 CNRS/Universit\'e de Nantes,
 2 rue de la Houssi\-ni\`e\-re,
 BP 92208,
 44322 Nantes Cedex 03, France.
}
\email{massuyea@math.univ-nantes.fr}

\begin{abstract}
We investigate classification results for general quadratic
functions on torsion abelian groups. Unlike the previously studied
situations, general quadratic functions are allowed to be
inhomogeneous or degenerate.  We study the discriminant construction
which assigns, to an integral lattice with a distinguished
characteristic form, a quadratic function on a torsion group. When
the associated symmetric bilinear pairing is fixed, we construct an
affine embedding of a quotient of the set of characteristic forms
into the set of all quadratic functions and determine explicitly its
cokernel.  We determine a suitable class of torsion groups so that
quadratic functions defined on them are classified by the stable
class of their lift. This refines results due to A.H. Durfee, V.
Nikulin, C.T.C. Wall and E. Looijenga -- J. Wahl. Finally, we show that on this class of torsion groups, two quadratic functions $q, q'$ are
isomorphic if and only if they have equal associated Gauss sums and
there is an isomorphism between the associated symmetric bilinear
pairings $b_{q}$ and $b_{q'}$ which sends $d_{q}$ to $d_{q'}$, where
$d_{q}$ is the homomorphism defined by $d_{q}(x) = q(x) - q(-x)$.
This generalizes a classical result due to V. Nikulin.
Our results are elementary in nature and motivated by low-dimensional topology.
\end{abstract}
\maketitle

 A quadratic function $q$ on an abelian group $G$ is a map, with
values in an abelian group, such that the map $b: (x,y) \mapsto
q(x+y)-q(x)-q(y)$ is $\Z$-bilinear. Such a map $q$ satisfies $q(0) =
0$. If, in addition, $q$ satisfies the relation $q(nx) = n^{2}q(x)$
for all $n \in \Z$ and $x \in G$, then $q$ is homogeneous. In
general, a quadratic function cannot be recovered from the
associated bilinear pairing $b$. Homogeneous quadratic functions on
torsion groups first appeared as quadratic enhancements of the
linking pairing on the torsion subgroup of the $(2n-1)$-th homology
group of an oriented $(4n-1)$-manifold. Typically, these quadratic
enhancements appear in topology when the manifold is equipped with a
framing \cite{BM} \cite{MS} \cite{LL}. They were used as a
fundamental ingredient in the classification up to regular homotopy
of immersed surfaces in $\R^{3}$ \cite{Pin}. They were extensively
studied from the algebraic viewpoint of Witt and Grothendieck
groups, see for instance
 \cite{Durfee} \cite{Ka} \cite{Lan}. However, there are topological
 motivations to consider inhomogeneous enhancements of the linking
 pairing \cite{LW} \cite{De2}. It is also convenient to consider possibly
 degenerate quadratic functions \cite{De1}. The motivation for considering
{\it{general}} quadratic functions stems from our work on closed
Spin$^c$-manifolds of dimension 3 and their finite type invariants
\cite{DM}.

This paper studies, and gives classification results for, quadratic
functions on torsion abelian groups with values in $\Q/\Z$. Those
results are relatively well known in the case of
{\emph{nondegenerate}} symmetric bilinear pairings and
{\emph{homogeneous}} quadratic functions. However, the authors have
not succeeded in finding in the litterature the general results
for quadratic functions. \\

To describe the first result, we review (\S
\ref{subsec:discriminant}) a construction, known as the discriminant
construction. This construction assigns to a symmetric bilinear
lattice $M$  equipped with a certain element (called a
characteristic form) $c \in {\rm{Hom}}(M,\Z)$, a quadratic function
on a torsion abelian group, with values in $\Q/\Z$. The first result
is an embedding theorem (Theorem \ref{th:algebraic-embedding}) which
describes which quadratic functions on a torsion abelian group with
values in $\Q/\Z$, come from a characteristic form (for a fixed
symmetric bilinear pairing). Unlike the previously studied
situations, the quadratic functions obtained by the discriminant
construction may be degenerate and the torsion abelian group may be infinite. \\

The second result addresses the problem of classifying up to
isomorphism the quadratic functions arising from the discriminant
construction. This problem is closely related to the stable
classification of symmetric bilinear lattices equipped with
characteristic forms. More precisely, we prove (\S
\ref{sec:stable}) that the stable equivalence on lattices is
equivalent to a particular class of isomorphism between the
associated quadratic functions, which we determine (Theorem
\ref{th:stabilisation}). The main application of this result has a
particularly simple form: the stable equivalence can be realized
by adding $1$-dimensional unimodular lattices (Corollary
\ref{cor:stabsimple}).
 This generalizes a classical result
for lattices without distinguished characteristic forms.\\

While the second result theoretically solves the classication
problem of a vast class of quadratic functions (which includes all
nondegenerate quadratic functions when the group is finite), or at
least reduces it to a stable classification problem for lattices
with characteristic forms, it does not seem to allow to decide
concretely whether two quadratic functions are isomorphic. This
leads to the problem of classification of quadratic functions by
means of invariants. We propose a solution to the problem on finite
abelian groups by giving a short proof
 (\S \ref{sec:system}, Theorem \ref{th:invcomplet}) that the
isomorphism class of a quadratic function $q$ on a finite abelian
group $G$ is determined by the isomorphism class of the associated
bilinear pairing $b$ equipped with a distinguished element $d_{q}
\in {\rm{Hom}}(G,\Q/\Z)$, called the homogeneity defect of $q$, and
the Gauss sum $\sum_{x \in G} \exp(2\pi iq(x))$ associated to $q$.
This result is then extended to a more general class of torsion
groups. The particular case of homogeneous quadratic functions
($d_{q} = 0$) on a finite group $G$ was proved by V. Nikulin by an
induction on the number of generators of $G$
\cite[Theorem 1.11.3]{Nik}. \\

These three results were applied to the classification
of degree 0 finite type invariants of Spin$^c$ $3$-manifolds
\cite{DM}.\\

\noindent{\bf{Acknowledgements}}. F.D. was supported by a E.U.
Marie Curie Research Fellowship (HMPH-2001-01174).

\tableofcontents

\section{Preliminaries} \label{sec:quad}

Let $G$ be an abelian group. A set $X$ on which $G$ acts freely and
transitively is an {\it{affine space over}} $G$. For such actions,
the left multiplicative notation is used. Any bilinear pairing $b: G
\times G' \to H$, where $G,\ G'$ and $H$ are abelian groups, has a
left (resp. right) adjoint map
$\hphantom{a}_{\centerdot}\widehat{b}:G \to \hbox{Hom}(G',H)$ (resp.
$\widehat{b}_{\centerdot}:G' \to \hbox{Hom}(G,H)$), defined by
$\hphantom{a}_{\centerdot}\widehat{b}(x)(y) = b(x,y)$ (resp.
$\widehat{b}_{\centerdot}(y)(x) = b(x,y)$), $x \in G,$ $y \in G'$.
We say that the bilinear pairing $b$ is {\it{left nondegenerate}}
(resp. {\it{left nonsingular}} or {\it{left unimodular}}) if its
left adjoint map is injective (resp. bijective). We define similarly
right nondegenerate and right nonsigular. The bilinear pairing $b$
is {\it{nondegenerate}} (resp. {\it{nonsingular}}) if $b$ is both
left and right nondegenerate. In the case when $G = G'$ and $b$ is
symmetric, $\hphantom{a}_{\centerdot}\widehat{b} =
\widehat{b}_{\centerdot}$ is denoted by $\widehat{b}$. Throughout
the paper, the abelian group $G$ will be a finitely generated free
module over $R = \Z$ or $\Q$, or a torsion group. We set $C_G = R$
if $G$ is free and $C_G = \Q/\Z$ if $G$ is torsion. It is convenient
to define the dual of $G$ by $G^{*} = \hbox{Hom}(G,C_G)$.

A map $q:G \to \Q/\Z$ on a torsion abelian group $G$ is a
{\it{quadratic function}} if that the associated map $b_{q}:G \times
G \to \Q/\Z$, defined by $b_{q}(x,y) = q(x+y) - q(x) - q(y)$, is a
bilinear pairing. Note that this definition implies that $q(0) = 0$
for any quadratic function $q$. A quadratic function $q:G \to \Q/\Z$
satisfying $q(nx) = n^{2}q(x)$ for all $x \in G$ and all $n \in \Z$
is said to be {\it{homogeneous}}. We say that a map $q:G \to \Q/\Z$
is {\it{quadratic over}} a symmetric bilinear pairing $b:G \times G
\to \Q/\Z$ if $b_{q} = b$. A quadratic function $q$ is
{\it{nondegenerate}} (resp. {\it{nonsingular}} or {\it{unimodular}})
if $b_{q}$ is. The difference between two quadratic functions $q,
q'$ over the same bilinear pairing lies in $G^{*}$. In particular,
the set ${\rm{Quad}}(b)$ of all quadratic functions $q: G \to \Q/\Z$
over a symmetric bilinear pairing $b:G \times G \to \Q/\Z$ is affine
over the group $G^{*}$. The group $G$ acts also on ${\rm{Quad}}(b)$
via the adjoint map $\widehat{b}: G \to G^{*}$: $$ G \times
{\rm{Quad}}(b) \to {\rm{Quad}}(b), \ \  (\alpha, q) \mapsto \alpha
\cdot q = q + \widehat{b}(\alpha).$$ If $q$ is nonsingular, then
${\rm{Quad}}(b)$ is affine over the group $G$. To any $q \in
{\hbox{Quad}}(b)$, we associate the quadratic functions $-q \in
{\hbox{Quad}}(-b)$ and $\bar{q} \in {\hbox{Quad}}(b)$ respectively
by $(-q)(x) = -q(x)$ and $\bar{q}(x) = q(-x)$, $x \in G$.  Clearly
the maps $q \mapsto -q$ and $q \mapsto \bar{q}$ are involutive
bijections. The {\it{homogeneity defect}} $d_{q}$ of a quadratic
function $q \in {\rm{Quad}}(b)$ is defined by $d_{q} = q - {\bar{q}}
\in G^{*}$. Since $$ n^{2}q(x)-q(nx) = \frac{n(n-1)}{2} \ d_q(x),$$
a quadratic function $q$ is homogeneous if and only if $d_{q} = 0$.

Let $q:G \to \Q/\Z$ be a quadratic function with $b_{q} = b$. Let
$\psi:G' \to G$ be a group isomorphism. We define $\psi^{*}b:
G'\times G' \to \Q/\Z$ to be the symmetric bilinear pairing
defined by
$$\psi^{*}b(x,y) = b(\psi(x),\psi(y))\ \
{\rm{for}}\ {\rm{all}}\ x, y \in G'.$$ Similarly we define the
quadratic function $\psi^{*}q:G'\to \Q/\Z$ by $$\psi^{*}q(x) =
q(\psi(x))\ \ {\rm{for}}\  {\rm{all}}\ x \in G'.$$ We say that two
symmetric bilinear pairings $b,\ b'$ (resp. two quadratic functions
$q$, $q'$) defined on $G$ and $G'$ with values in $\Q/\Z$ are
{\it{isomorphic}}, and we write $b \sim b'$ (resp. $q \sim q'$), if
there exists a group isomorphism $\psi:G' \to G$ such that
$\psi^{*}b = b'$ (resp. $\psi^{*}q = q'$). An isomorphism
$\psi:G'\to G$ induces a bijective correspondence between
${\hbox{Quad}}(b)$ and ${\hbox{Quad}}(\psi^{*}b)$. In particular,
the subgroup Iso$(b)$ of automorphisms of $G$ preserving $b$ acts on
${\hbox{Quad}}(b)$: $\psi \cdot q = \psi^{*}q$.

\newpage
\section{Presentations of quadratic functions}
\label{sec:presentations}
\subsection{Lattices and their characteristic forms}
 \label{subsec:latticesandforms}
A {\emph{lattice}} $M$ is a finitely generated free abelian group. A
\emph{(symmetric) bilinear lattice} $(M,f)$ is a symmetric bilinear
form $f: M\times M \to \Z$ on a lattice $M$. Consider the vector
space  $V = M \otimes \Q$ over $\Q$. The dimension of $V$ is finite
and equal to the rank of $M$. Any bilinear lattice $(M,f)$ gives
rise by extension of scalars to a bilinear pairing $f_{\Qi}:V \times
V \to \Q$. Let $M^\sharp=\{ x\in V : f_{\Qi}(x,M) \subset \Z\}$ be
the \emph{dual lattice} for $(M,f)$. Clearly $M \subset M^{\sharp}$.
%
%
More generally, for any subgroup $N$ of $V$, we can define
$N^{\sharp} =\{ x\in V:  f_{\Qi}(x,N) \subset \Z\}$. A
\emph{fractional} (resp. \emph{integral}) \emph{Wu class} for
$(M,f)$ is an element $w \in V$ (resp. an element $w \in M$) such
that
\begin{displaymath}
\forall x\in M,\  f(x,x)-f_{\Qi}(w,x) \in 2\Z.
\end{displaymath}
The set of fractional (resp. integral) Wu classes is denoted by
Wu$^{\Qi}(f)$ (resp. Wu$(f)$) and is contained in $M^\sharp$.
Furthermore, Wu$^{\Qi}(f)$ is an affine space over $2M^{\sharp}$.\\
A \emph{characteristic form} for $f$ is an element $c \in M^*$
satisfying
\begin{displaymath}
\forall x\in M,\ f(x,x)-c(x) \in 2\Z. \label{eq:chardef}
\end{displaymath}
The set Char$(f)$ of characteristic forms for $f$ is an affine space
over $\hbox{Hom}(M,2\Z)$. Since $x \mapsto f(x,x)$ mod $2\Z$ is a
homomorphism, Char$(f)$ is not empty. Each fractional Wu class gives
a characteristic form by the equivariant map $w \mapsto
f_{\Qi}(w,-)|_{M}$, ${\hbox{Wu}}^{\Qi}(f) \to {\rm{Char}}(f)$.
\begin{lemma}
\label{lem:nondegenerate} If $f$ is nondegenerate, then the map $w
\mapsto f_{\Qi}(w,-)|_{M}$ is a bijective correspondence between
{\rm{Wu}}$^{\Qi}(f)$ and ${\rm{Char}}(f)$.
\end{lemma}
\begin{proof}
Since $f$ is nondegenerate, $\widehat{f}:M \to M^{*}$ is injective.
Then $\widehat{f}_{\Qi}: V \to V^{*} = {\rm{Hom}}(V,\Q)$ is
bijective. Hence, the map  $x \mapsto f_{\Qi}(x,-)|_{M},\
{\hbox{Wu}}^{\Qi}(f) \to {\hbox{Char}}(f)$ is also bijective.
\end{proof}
The quotient $\bar{M} = M/\hbox{Ker}\ \widehat{f}$ is finitely
generated free. Hence, the short exact sequence
\begin{equation}
\label{eq:splitting_of_M}
\xymatrix{
0 \ar[r] &  {\hbox{Ker}}\ {\widehat{f}} \ar[r]^-{i} & M \ar[r]^-{p} &
{\bar{M}} \ar[r] & 0}
\end{equation}
is split. Any section $s$ of $p$ induces an isomorphism
$$
(M,f) \simeq \left(\bar{M},\bar{f}\right) \oplus
\left(\hbox{Ker}\ \widehat{f},0\right),
$$
where $\bar{f}: \bar{M} \times \bar{M} \to \Z$ is the nondegenerate
pairing induced by $f$.
\begin{lemma}
\label{lem:retract} There is an injection $p^{*}|_{{\rm{\scriptsize
Char}}\left(\bar{f}\right)}: {\rm{Char}}\left(\bar{f}\right) \to
{\rm{Char}}(f)$ induced by $p$, and any section $s$ of $p$ induces
an affine retraction $s^{*}|_{{\rm{\scriptsize
Char}}(f)}:{\rm{Char}}(f) \to {\rm{Char}}(\bar{f})$ for
$p^{*}|_{{\rm{\scriptsize Char}} (\overline{f})}$.
\end{lemma}
\begin{proof}
Since $p\circ s=\hbox{Id}_{\bar{M}}$, we have $s^* \circ
p^*=\hbox{Id}_{\bar{M}^*}$. Moreover, one easily verifies that
$p^*\left(\hbox{Char}\left(\bar{f}\right)\right)
\subset\hbox{Char}\left(f\right)$ and
$s^*\left(\hbox{Char}\left(f\right)\right)\subset
\hbox{Char}\left(\bar{f}\right)$.
\end{proof}

\subsection{The discriminant construction}
\label{subsec:discriminant}

Suppose we are given a bilinear lattice $(M,f)$ as above. Consider
the torsion group $G_f = M^{\sharp}/M$ and the symmetric bilinear
pairing
$$L_f: G_f \times G_f \to \Q/\Z$$
defined by
\begin{equation}
\label{eq:pairing} L_f\left([x],[y]\right) = f_{\Qi}(x,y)\
\hbox{mod}\ \Z, \ \ x, y \in M^{\sharp}.
\end{equation}
Since $M^{\sharp\sharp}=M+ \hbox{Ker}\ \widehat{f}_{\Qi}$, the
radical of $L_f$ is
\begin{equation}
\hbox{Ker}\ \widehat{L}_f = \left(M+\hbox{Ker}\
\widehat{f}_{\Qi}\right)/M =\hbox{Ker}\ \widehat{f}_{\Qi}/
\hbox{Ker}\ \widehat{f} = \left(\hbox{Ker}\ \widehat{f}\right)
\otimes \Q/\Z. \label{eq:radical}
\end{equation} In particular, $L_f$ is nondegenerate if and only if $f$ is nondegenerate.\\

Consider the torsion subgroup $\hbox{T}\ {\hbox{Coker}}\
\widehat{f}$ of Coker $\widehat{f}$. The adjoint map
$\widehat{f}_{\Qi}:V \to V^{*}$ restricted to $M^{\sharp}$ induces
a canonical epimorphism $B_{f}:G_f \to \hbox{T}\ {\hbox{Coker}}\
\widehat{f}$. Hence there is a short exact sequence
\begin{equation}
0 \to  \hbox{Ker}\ \widehat{L}_f \to G_f \to \hbox{T}\
{\hbox{Coker}}\ \widehat{f} \to 0. \label{eq:shortexact}
\end{equation}
It follows that the symmetric bilinear form $L_f$ factors
to a nondegenerate pairing
$$\lambda_f: {\hbox{T}}\ {\hbox{Coker}}\ \widehat{f}
\times {\hbox{T}}\ {\hbox{Coker}}\ \widehat{f} \to \Q/\Z.$$ It
also follows from (\ref{eq:shortexact})  that $G_{f}=0$ if and
only if $M^{\sharp}=M$ if and only if $f$ is unimodular. The map
$B_{f}:G_{f} \to {\rm{T}}\ {\rm{Coker}}\ \widehat{f}$ induces a
canonical isomorphism $\Bigl({\rm{T}}\ {\rm{Coker}}\ \widehat{f},
\lambda_{f}\Bigr) \simeq \Bigl(G_{\bar{f}}, L_{\bar{f}}\Bigr)$. In
general, observe that a section $s$ of $p$ in
$(\ref{eq:splitting_of_M})$ induces a section of
(\ref{eq:shortexact}). Hence there is a (non-canonical) orthogonal
decomposition:
 \begin{equation}
 (G_f, L_f) \simeq \left({\hbox{T}}\ {\hbox{Coker}}\ \widehat{f},
 \lambda_f\right) \oplus \left(\hbox{Ker}\ \widehat{L}_f, 0\right).
 \label{eq:splitq}
 \end{equation}
For future use, we notice that the converse of our previous observation holds.

\begin{lemma} \label{lem:sections}
Any section of $B_{f}:G_{f} \to {\rm{T}}\ {\rm{Coker}}\ \widehat{f}$ is induced by
a section of $p:M \to \bar{M}$.
\end{lemma}

\begin{proof}
The underlying map of the affine map from sections of $p$ to
section of $B_{f}$ is the natural homomorphism
${\rm{Hom}}(\bar{M},{\rm{Ker}}\ \widehat{f}) \to
{\rm{Hom}}(G_{\bar{f}}, {\rm{Ker}}\ \widehat{L}_{f})$. It suffices
to lift $\alpha \in {\rm{Hom}}(G_{\bar{f}}, {\rm{Ker}}\
\widehat{L}_{f})$ to a map $\tilde{\alpha} \in
{\rm{Hom}}(\bar{M}^{\sharp},{\rm{Ker}}\ \widehat{f}_{\Qi})$; the
restriction $\tilde{\alpha}|_{\bar{M}}$ will be the desired
section of $p$. Since ${\bar{M}}^{\sharp}$ is free, the lift
exists.
\end{proof}

Suppose next that $(M,f,c)$ is a bilinear lattice equipped with a
characteristic form $c \in M^{*}$. Denote by $c_{\Qi}:V\to \Q$ the
linear extension of $c$ . We associate to $(M,f,c)$ a quadratic
function $\phi_{f,c}: G_f \to \Q/\Z$ over $L_f$ by
\begin{displaymath}
\phi_{f,c}\left([x]\right)=\frac{1}{2}(f_{\Qi}(x,x)-c_{\Qi}(x))\
\hbox{mod}\ \Z, \ \ x \in M^{\sharp}.
\end{displaymath}
Observe that the homogeneity defect $d_{\phi_{f,c}} \in
G_{f}^{*}=\hbox{Hom}\left(G_f,\Q/\Z\right)$ is given by
\begin{equation}
d_{\phi_{f,c}}\left([x]\right) = -c_{\Qi}(x)\ \hbox{mod}\ \Z, \ \
x \in M^{\sharp}. \label{eq:homogeneity-defect}
\end{equation}
\begin{definition}
The triple $(M,f,c)$ is said to be a \emph{presentation} of the
quadratic function $\phi_{f,c}$ on the torsion group $G_f$. The
assignation $(M,f,c)\mapsto (G_f,\phi_{f,c})$ is called the
\emph{discriminant} construction.
\end{definition}
\begin{lemma} \label{lem:orthogonalsum}
The discriminant construction preserves orthogonal sums.
\end{lemma}

\begin{remark}
\label{rem:particular_case_known}
The following conditions are seen to be equivalent from
(\ref{eq:radical}) and (\ref{eq:shortexact}):
\begin{enumerate}
\item[$\cdot$] $f$ is nondegenerate;
\item[$\cdot$] $L_{f}$ is nondegenerate;
\item[$\cdot$] $G_{f}$ is a finite product of cyclic groups;
\item[$\cdot$] Coker $\widehat{f}$ is a finite product of cyclic groups;
\item[$\cdot$] $G_{f}$ and Coker $\widehat{f}$ are isomorphic.
\end{enumerate}
If one of those conditions is satisfied, $B_{f}: G_f \to
\hbox{Coker}\ \widehat{f}$ is an isomorphism which induces a
bijective correspondence $\hbox{Quad}\left(\lambda_f\right)\simeq
\hbox{Quad}\left(L_f\right).$ Hence $\phi_{f,c}$ factors via
$B_{f}$ to a nondegenerate quadratic function over $\lambda_{f}$.
More generally, $\phi_{f,c}$ factors via $B_{f}$ to a
nondegenerate quadratic function over $\lambda_{f}$ if and only if
 $c$ can be taken to be the image of a fractional Wu class
(cf. Lemma \ref{lem:nondegenerate}). This particular case coincides with the
 usual discriminant construction, as presented for example in
\cite[\S 1.2.1]{De1}.

The map $(M,f,c) \mapsto (G_f, \phi_{f,c})$, from the monoid of
nondegenerate bilinear lattices equipped with characteristic forms
to the monoid of isomorphism classes of nondegenerate quadratic
functions on finite groups, is known to be surjective \cite{Wal1}.
\end{remark}

\begin{remark}
There exists a discriminant construction which applies to
arbitrary Dedekind rings (but produces only {\it{homogeneous}}
quadratic functions) instead of $\Z$, see \cite{Durfee}. This
leads to the following\\

\noindent{\bf{Question}}. Does there exists a generalization of the
discriminant construction producing from lattices over Dedekind
domains quadratic functions on torsion modules over Dedekind domains
? (The notion of quadratic function should be generalized first.)
More generally, can one define a localization of quadratic functions
? (For homogeneous quadratic functions, see \cite{Ka}.)
\end{remark}

\subsection{Properties of the discriminant construction}
Consider the bilinear map  $M^{*} \times M^{\sharp} \to \Q$
defined by $(\alpha, x) \mapsto \alpha_{F}(x)$. This map induces a
bilinear pairing
$$\langle -, -\rangle: \ \hbox{Coker}\ \widehat{f} \times G_{f} \to \Q/\Z.$$
\begin{lemma} \label{lem:kernels}
The bilinear pairing $\langle - , - \rangle: {\rm{Coker}}\
\widehat{f} \times G_{f} \to \Q/\Z$ is left nondegenerate
 (respectively left nonsingular if and only if $f$ is nondegenerate)
and right nonsingular.
\end{lemma}

\begin{proof}
Consider the left adjoint map Coker $\widehat{f} \to
\hbox{Hom}(G_{f}, \Q/\Z), [\alpha] \mapsto \langle [\alpha], -
\rangle$ (mod $\Z$). Let $[\alpha]$ lie in the kernel, that is,
$\alpha_{\Qi}\left(M^{\sharp}\right) \subseteq \Z$. To prove that
$[\alpha] = 0$, it is sufficient to show that  there exists $v \in
M$ such that $\alpha=\widehat{f}(v)$. Assume first that $f$ is
nondegenerate.
 Then $\widehat{f}_{\Qi}$ is bijective,
so there is an element $v \in V$ such that $\alpha_{\Qi} =
\widehat{f}_{\Qi}(v)$. Since $\alpha_{\Qi}(M^{\sharp}) =
f_{\Qi}(v, M^{\sharp}) \subset \Z$, this means that $v \in
M^{\sharp \sharp} = M$. To deal with the general case ($f$
possibly degenerate), we consider
 the nondegenerate symmetric bilinear
pairing $\bar{f}$ induced by $f$ on the lattice
$\bar{M} = M/\hbox{Ker}\ \widehat{f}$. We shall use the
following observation.

\begin{sublemma}
Let $N$ be a free abelian group and let $W = N \otimes \Q$. Let
$\gamma \in {\rm{Hom}}_{\Qi}(W,\Q)$. We have $\gamma = 0$ if and
only if $\gamma(W) \subseteq \Z$.
\end{sublemma}

Recall that Ker $\widehat{f}$ is a free subgroup of $M$ and Ker
$\widehat{f}_{\Qi} = ({\rm{Ker}}\ \widehat{f}) \otimes \Q$. We
have $\alpha_{\Qi}({\rm{Ker}}\ \widehat{f}_{\Qi}) \subseteq
\alpha_{\Qi}(M^{\sharp}) \subseteq \Z$. Apply the sublemma above
with $N = {\rm{Ker}}\ \widehat{f}$ to $\alpha_{\Qi}|_{ {\rm{Ker}}\
\widehat{f}_{\Qi}}$ yields $\alpha_{\Qi}({\rm{Ker}}\
\widehat{f}_{\Qi}) = 0$. In particular, $\alpha( {\rm{Ker}}\
\widehat{f}) = 0$. Hence, $\alpha$ induces a form
$\bar{\alpha}:\bar{M} \to \Z$ such that
$\bar{\alpha}\left(\bar{M}^\sharp\right) \subseteq \Z$. Applying
the previous (nondegenerate) case to $\bar{\alpha}$ yields an
element $\bar{v} \in \bar{M}$ such that $\bar{\alpha} =
\widehat{\bar{f}}\left(\bar{v}\right)$. Any lift $v \in M$ of
$\bar{v}$ satisfies $\alpha = \widehat{f}(v)$. So $[\alpha] = 0$
as desired.

Suppose that $f$ is nondegenerate. By Remark
$\ref{rem:particular_case_known}, G_{f}$ (resp.
Coker~$\widehat{f}$) is a product of cyclic groups
$\frac{1}{n}\Z/\Z$ (resp. $\Z/n\Z$). The injective left adjoint
map Coker $\widehat{f} \to \hbox{Hom}(G_{f}, \Q/\Z)$ is between
two groups of the same finite order, so must be surjective as
well. Suppose that $f$ is not nondegenerate. Then
Coker~$\widehat{f}$ (resp. $G_{f}$) has a summand $\Z$ (resp. a
summand which is $\Q/\Z$). It suffices to show that the map $\Z
\to \hbox{Hom}(\Q/\Z, \Q/\Z),\ r \mapsto r\ {\rm{Id}}_{\Qi/\Zi}$
is not surjective. Recall that $\Q/\Z = \oplus_{p} A_{p}$, where
the direct sum is over all primes and $A_p$ is the subgroup
consisting of elements $x \in \Q/\Z$ such that $p^{k}x = 0$ for
some $k$. Consider the map $\pi_{p}: A_{p} \to A_{p}$ which is
multiplication by $p$. Then the product $\prod_{p} \pi_{p} \in
\prod_{p} \hbox{Hom}(A_{p},\Q/\Z) = \hbox{Hom}(\Q/\Z,\Q/\Z)$ is
not in the image.

Consider the right adjoint map $G_{f} \to
\hbox{Hom}\left(\hbox{Coker}\ \widehat{f},\Q/\Z\right)$. An
element $[x] \in G_{f}$ lies in the kernel if and only if
$\alpha_{\Qi}(x) = 0$ mod $\Z$ for all $\alpha \in M^{*}$. In
particular, for a unimodular bilinear pairing $g$ on $M$,
$g_{\Qi}(v,x) = 0$ mod $\Z$ for all $v \in V$. Unimodularity of
$g$ implies that $x \in M$. Thus $[x] = 0$.

We now prove that the right adjoint map is surjective. Consider
the canonical isomorphism $V \to V^{**}, v \mapsto \langle -, v
\rangle$. This map sends $M^{\sharp}$ onto $N = \{ \langle -, w
\rangle \in V^{**}\ : \ \langle \widehat{f}(M), w \rangle
\subseteq \Z \}$ which is a subgroup of $V^{**}$. The subgroup
$N^{0} = \{ x \in \hbox{Hom}(M^{*},\Q)\ : \ \langle
\widehat{f}(M),x \rangle \subseteq \Z \}$ embeds in $N$ by linear
extension over $\Q$: $M^{*} \to M^{*} \otimes \Q = V^{*}$. Any $x
\in N^{0}$ induces a homomorphism $\hbox{Coker}\ \widehat{f} \to
\Q/\Z$, hence there is an induced homomorphism $N^{0} \to
\hbox{Hom}\left(\hbox{Coker}\ \widehat{f}, \Q/\Z\right)$, which is
obviously surjective (since $M^*$ is free). Hence any $\varphi \in
\hbox{Hom}\left(\hbox{Coker}\ \widehat{f}, \Q/\Z\right)$ lifts to
an element $x \in N^{0} \subseteq N$ where $x = \langle -, v
\rangle$ for some $v \in M^{\sharp}$. Thus $\varphi = \langle -,
[v] \rangle$.
\end{proof}
\begin{lemma} \label{lem:cokernel=ker*}
The cokernel of the inclusion ${\rm{T}}\ {\rm{Coker}}\ \widehat{f}
\hookrightarrow {\rm{Coker}}\ \widehat{f}$ is
$\left({\rm{Ker}}\ \widehat{f}\right)^{*}$.
\end{lemma}
\begin{proof}
Since $M/\hbox{Ker}\ \widehat{f}$ is free, the inclusion Ker
$\widehat{f} \hookrightarrow M$ induces a surjective homomorphism
$M^{*} \to \left(\hbox{Ker}\ \widehat{f}\right)^{*}$ which sends
$\widehat{f}(M)$ to $0$. Hence we obtain a well-defined surjective
homomorphism $\hbox{Coker}\ \widehat{f} \to \left(\hbox{Ker}\
\widehat{f}\right)^{*}$, the kernel of which, since
$\left(\hbox{Ker}\ \widehat{f}\right)^*$ is free of same rank as
$\hbox{Coker}\ \widehat{f}$, consists of all torsion elements.
\end{proof}
By $(\ref{eq:radical})$, any homomorphism Ker  $\widehat{f} \to
\Z$ induces by tensoring with $\Q/\Z$ a homomorphism $\hbox{Ker}\
\widehat{L}_f \to \Q/\Z$. Denote by $j_f:\left(\hbox{Ker}\
\widehat{f}\right)^{*} \to \left(\hbox{Ker}\
\widehat{L}_f\right)^{*}$ the corresponding homomorphism.

Let $c\in {\rm{Char}}(f)$. Observe that $\phi_{f,c+2\hat{f}(u)} =
\phi_{f,c}$ for any $u \in M$. Hence $\phi_{f,c}$ depends on $c$
only mod $2\widehat{f}(M)$. Furthermore, the abelian group
$M^{*}/\widehat{f}(M) = \hbox{Coker}\ \widehat{f}$ acts freely and
transitively on Char$(f)/2\widehat{f}(M)$ by
$$
[\alpha] \cdot [c] = [c + 2\alpha] \in \hbox{Char}(f)/2\widehat{f}(M), \
\alpha \in M^{*},\ c \in \hbox{Char}(f).
$$
The following result is the main result of this section.
\begin{theorem} \label{th:algebraic-embedding}
The map $c \mapsto \phi_{f,c}$ induces an affine embedding
$$\phi_{f}:{\rm{Char}}(f)/2\widehat{f}(M) \hookrightarrow
{\rm{Quad}}(L_{f})$$ over the opposite of the group monomorphism
$\hphantom{a}_{\centerdot}\langle -, - \rangle :{\rm{Coker}}\
\widehat{f} \hookrightarrow {\rm{Hom}}(G_{f},\Q/\Z)$. Furthermore,
$$ {\rm{Coker}}\ \phi_{f} \simeq {\rm{Coker}}\ j_{f} =
\frac{({\rm{Ker}}\ \widehat{L}_{f})^{*}}{j_{f}\Bigl(({\rm{Ker}}\
\widehat{f})^{*}\Bigr)} $$
and given $q \in {\rm{Quad}}(L_{f})$, the following two assertions are equivalent:
\begin{enumerate}
\item[$\cdot$] $q \in {\rm{Im}}\ \phi_{f}$;
\item[$\cdot$] $q|_{{\rm{\scriptsize Ker}}\ \widehat{L}_{f}} \in {\rm{Im}}\ j_{f}$.
\end{enumerate}
\end{theorem}

As a consequence of Lemma \ref{lem:kernels}, we have:

\begin{corollary} \label{cor:bijective-embedding}
The map $\phi_{f}$ is bijective if and only if $f$ is nondegenerate.
\end{corollary}

Combining Corollary \ref{cor:bijective-embedding} and Lemma \ref{lem:retract},
we obtain the following result.

\begin{corollary}
 For any section $s$ of $p:M \to \overline{M}$, the
 map $${\rm{Char}}(f) \to  {\rm{Quad}}(\lambda_{f}), \ c
\mapsto \phi_{\bar{f}, s^{*}|_{{\rm{Char}}(f)}(c)}$$ is onto.
\end{corollary}

\begin{remark}
This construction applies in particular when $(M,f)$ is the homology
group $H_2(X)$ of a compact connected oriented simply-connected
$4$-manifold $X$, equipped with its symmetric bilinear intersection
pairing. Then $\lambda_{f}$ can be identified with the torsion
linking pairing on T$H_{1}(\partial X)$ (up to sign depending on the
orientation of $(X, \partial X)$). For details, see for example
\cite{Durfee} \cite{De1}. Furthermore $G_{f}$ identifies with
$H_{2}(\partial X;\Q/\Z)$ and the affine space
${\hbox{Char}}(f)/2\widehat{f}(M)$ identifies with the affine space
of Spin$^c$-structures on $\partial X$. See \cite{DM} for further
details about the last point.
\end{remark}

\begin{proof}[Proof of Theorem $\ref{th:algebraic-embedding}$]
Let us prove that the map $\phi_{f}$ is indeed affine over the
homomorphism stated above. Let $\alpha \in M^{*}$ and $x \in
M^{\sharp}$: $\phi_{f,[\alpha] \cdot [c]}\left([x]\right) -
\phi_{f,[c]}\left([x]\right)=\phi_{f,c+2\alpha}\left([x]\right) -
\phi_{f,c}\left([x]\right) = -\alpha_{\Qi}(x)\ \hbox{mod}\ \Z$.
Hence $\phi_{f,[\alpha] \cdot [c]}([x])= \phi_{f,[c]}([x])-\langle
[\alpha],[x]\rangle.$ The fact that $\phi_{f}$ is injective follows
from the fact that the map $\hphantom{a}_{\centerdot}\langle -, -
\rangle:\hbox{Coker}\ \widehat{f} \to \hbox{Hom}\left(G_{f},
\Q/\Z\right)$ is injective and this was proved in Lemma
\ref{lem:kernels}.

Next, we proceed to determine Coker $\phi_{f}$. Consider the following diagram:
$$\xymatrix{
0 \ar[r] & {\hbox{T}}\ {{\hbox{Coker}}\ \widehat{f}} \ar[r]
\ar[d]^-{\widehat{\lambda}_{f}}_{\simeq} & {{\hbox{Coker}}\
\widehat{f}} \ar[r] \ar[d]^{\hphantom{a}_{\centerdot}\langle -, -
\rangle} & \left({{\hbox{Ker}}\
\widehat{f}}\right)^{*} \ar[r] \ar[d]^{j_{f}} & 0 \\
0 \ar[r] & \left({\hbox{T}}\ {{\hbox{Coker}}\
\widehat{f}}\right)^{*} \ar[r] & {\hbox{Hom}}(G_{f}, \Q/\Z) \ar[r]
& \left({\hbox{Ker}\ \widehat{L}_{f}}\right)^{*} \ar[r] & 0 }
$$
The first row is given by Lemma \ref{lem:cokernel=ker*}. The
second row is obtained from the short exact sequence
$(\ref{eq:shortexact})$ by applying the exact functor
Hom$(-,\Q/\Z)$ (since $\Q/\Z$ is divisible). It follows from the
definitions that the diagram is commutative.

A standard argument (for instance, the ``snake lemma'')
applied to the diagram leads to
Coker$\hphantom{a}_{\centerdot}\langle -, - \rangle \simeq \hbox{Coker}\ j_{f}$.

Finally we prove the last statement of the theorem. By definition,
$c($Ker $\widehat{f}) \subseteq 2\Z$. In particular,
$\phi_{f,c}|_{\hbox{\scriptsize Ker}\ \widehat{L}_{f}}:\hbox{Ker}\
\widehat{L}_{f} \to \Q/\Z$ is a homomorphism given by
\begin{equation}
\phi_{f,c}([x]) = - \frac{1}{2}c_{\Qi}(x)\ \hbox{mod}\ \Z, \ \ x \in
M + {\hbox{Ker}}\ \widehat{f}_{\Qi} \subseteq M^{\sharp}.
\label{eq:on-the-kernel}
\end{equation}
Hence $\phi_{f,c}|_{\hbox{\scriptsize Ker}\ \widehat{L}_{f}} =
j_{f}(-\frac{1}{2}c) \in {\hbox{Im}}\ j_{f}$. Conversely, let $q
\in \hbox{Quad}(L_{f})$ such that $q|_{\hbox{\scriptsize Ker}\
\widehat{L}_f} \in {\hbox{Im}}\ j_{f}$. Pick $c \in
\hbox{Char}(f)$.  Then the restriction homomorphism $(G_{f})^{*}
\to (\hbox{Ker}\ \widehat{L}_{f})^{*}$ induced by
Ker $\widehat{L}_{f} \hookrightarrow G_{f}$ sends $q - \phi_{f,c}$
into Im $j_{f}$. Since
$$\hbox{Coker}\hphantom{a}_{\centerdot}\langle -, - \rangle
\simeq \hbox{Coker}\ j_{f},$$ there is an element $[\alpha] \in
\hbox{Coker}\ \widehat{f}$ such that $q - \phi_{f,c} = \langle
[\alpha], - \rangle$. Since the map $\phi_{f}$ is affine over
$\hbox{Coker}\ \widehat{f}\hookrightarrow \hbox{Hom}\left(G_{f},
\Q/\Z\right)$, it follows that $q = \phi_{f, [-\alpha] \cdot
[c]}.$
\end{proof}

\begin{corollary} \label{co:which-homomorphisms}
A homomorphism $h:G_{f} \to \Q/\Z$ is in the image of
$\hphantom{a}_{\centerdot}\langle -, - \rangle : {\rm{Coker}}\
\widehat{f} \to {\rm{Hom}}(G_{f},\Q/\Z)$ if and only if
$h|_{{\rm{Ker}}\ \widehat{L}_{f}}:{\rm{Ker}}\ \widehat{L}_{f} \to
\Q/\Z$ lifts to a homomorphism ${\rm{Ker}}\ \widehat{f} \to \Z$.
\end{corollary}

\begin{proof}
Set $h'= h|_{{\rm{Ker}}\ \widehat{L}_{f}}$.  By definition,  $h'$
lifts to Ker $\widehat{f}$ if and only if $h' \in {\rm{Im}}\ j_{f}$.
Now apply Theorem \ref{th:algebraic-embedding}.
\end{proof}

We complete this section with a simple observation on Coker
$\phi_{f}$. Let $p$ be a prime. Let $\displaystyle \Z_p =
\lim_{\leftarrow} \Z/p^{k}\Z$ be the ring of $p$-adic numbers. Set
$ \displaystyle \widehat{\Z} = \prod_{p}
 \Z_{p}$ where the product is taken over all primes. The map $n \mapsto (n \
 {\rm{mod}}\ p^{k})_{k \geq 0}$ defines a
 natural embedding $\Z \hookrightarrow {\Z}_{p}$.
 Consequently,  there is a diagonal embedding $\Z \hookrightarrow \widehat{\Z}$.
 Since $\Q/\Z$ is a direct
limit of finitely generated subgroups and since any finitely
generated torsion group decomposes as a direct sum of cyclic
subgroups, we see the following

\begin{lemma} \label{eq:comp}
There is an isomorphism $\psi:{\rm{Hom}}(\Q/\Z,\Q/\Z) \simeq
\widehat{\Z}.$
\end{lemma}

Denote by $j$ the embedding $\Z \to {\rm{Hom}}(\Q/\Z,\Q/\Z), \ 1
\mapsto {\rm{Id}}_{\Qi/\Zi}$. The diagram
$$ \xymatrix{
{\Z} \ar[rd]^-{j} \ar[d] & \\
{\widehat{\Z}} \ar[r]_-{\psi} & {\rm{Hom}}(\Q/\Z,\Q/\Z)\\
 }$$
is commutative. Recall, from $(\ref{eq:radical})$, that $j_{f}$ is
the canonical injection ${\rm{Hom}}({\rm{Ker}}\ \widehat{f}, \Z) \to
\linebreak {\rm{Hom}}({\rm{Ker}}\ \widehat{f} \otimes \Q/\Z,
\Q/\Z)$. Since Ker $\widehat{f}$ is free abelian, we may regard (by
means of Lemma \ref{eq:comp}) $j_{f}$ simply as the diagonal
embedding $\delta : {\rm{Ker}}\ \widehat{f} \to {\rm{Ker}}\
\widehat{f} \otimes \widehat{\Z}$. Therefore,

\begin{equation}
{\rm{Coker}}\ \phi_{f} \simeq {\rm{Coker}}\ j_{f} \simeq  \frac{\Bigl({\rm{Ker}}\
{\widehat{f}}\Bigr) \otimes \widehat{\Z}}{\delta\Bigl({\rm{Ker}}\
{\widehat{f}}\Bigr)}.
\end{equation}

\section{The stable classification theorem}
\label{sec:stable}

The goal is to generalize a result due to Wall and Durfee
(\cite[Corollary 1]{Wal2} and \cite[Theorem 4.1]{Durfee}).
We define a natural notion of stable equivalence on lattices equipped
with a characteristic form. The resulting stable classification problem
is shown to be essentially equivalent to the classification of
the quadratic functions
induced by the discriminant construction (\S \ref{sec:presentations}). \\

There is a natural notion of isomorphism among triples $(M,f,c)$
defined by bilinear lattices with characteristic forms: we say that two
triples $(M,f,c)$ and $(M',f',c')$ are isomorphic (denoted $(M,f,c)
\simeq (M',f',c'))$ if there is an isomorphism
$\psi:M \to M'$ such that $\psi^{*}f' = f$ and $\psi^{*}c' = c$ mod
$2\widehat{f}(M)$.  All such triples form a monoid for the orthogonal sum
$\oplus$. Two triples $(M,f,c)$ and $(M',f',c')$ are said to be
\emph{stably equivalent} if they become isomorphic after stabilizations
with some unimodular lattices, that is, there is an isomorphism
(called a {\emph{stable equivalence}}) between $(M,f,c) \oplus (U,g,u)$ and
$(M',f',c') \oplus (U',g',u')$ for some unimodular lattices $(U,g,u)$ and
$(U',g',u')$ equipped with characteristic forms $u$ and $u'$ respectively.
Since, as we have seen in \S \ref{sec:presentations}, unimodular lattices
induce trivial discriminant bilinear forms and quadratic functions,
a stable equivalence between two triples induces an isomorphism of the
corresponding quadratic functions. We are interested in whether the converse
holds and to which extent. In fact, a positive answer is provided
in the case of nondegenerate lattices.
%
%
%
\begin{proposition} \label{prop:stabilisation-nondegeneree}
Two nondegenerate symmetric bilinear lattices $(M,f,c)$ and $(M',f',c')$
equipped with characteristic forms are stably equivalent if and only if
their associated quadratic functions $(G_{f}, \phi_{f,c})$ and
$(G_{f'}, \phi_{f',c'})$ are isomorphic. In fact,
any isomorphism $\psi$ between their associated quadratic functions
$(G_{f}, \phi_{f,c})$ and $(G_{f'}, \phi_{f',c'})$ lifts
to a stable equivalence between  $(M,f,c)$ and $(M',f',c')$.
\end{proposition}

\begin{proof}
It is sufficient to prove the second statement. There is an
obvious notion of stable equivalence for symmetric bilinear
lattices (simply forget the characteristic form). Consider the
symmetric bilinear lattices $(M,f)$ and $(M',f')$. The hypothesis
implies that $L_{f'} = \psi^{*}L_{f}$. It follows from \cite[Proof
of Th. 4.1]{Durfee} that there is stable equivalence ${\mathfrak{s}}$ between
$(M,f)$ and $(M',f')$ inducing $\psi$. Explicitly, ${\mathfrak{s}}$ is an
isomorphism $(M,f) \oplus (U,g) \simeq (M',f') \oplus (U',g')$,
where $(U,g)$ and $(U',g')$ are unimodular lattices. In
particular, ${\mathfrak{s}}$ induces an affine isomorphism Char$(f \oplus g)
\simeq {\rm{Char}}(f' \oplus g')$. Let $u \in \hbox{Char}(g)$. We have
 ${\mathfrak{s}}(c \oplus u) =
c''
\oplus u' \in {\rm{Char}}(f' \oplus g') = {\rm{Char}}(f') \times
{\rm{Char}}(g')$. We have to show
that $c'' = c' \ {\rm{mod}}\ 2\widehat{f'}(M')$. By
 hypothesis, $\phi_{f',c'} \circ \psi =  \phi_{f,c}$. Since
$\psi$ is induced by ${\mathfrak{s}}$, we have
$\phi_{f',c'} = \phi_{f',{\mathfrak{s}}(c)}$. The
desired result follows from Theorem \ref{th:algebraic-embedding}.
\end{proof}

In the general case, we have to deal with the potential degeneracy
of lattices. The first observation is that it is {\it{not true}}
that any isomorphism between $(G_{f}, \phi_{f,c})$ and $(G_{f'},
\phi_{f',c'})$ will lift to a stable equivalence between $(M,f,c)$
and $(M',f',c')$. In fact, the simplest counterexample is given by
the degenerate bilinear lattice $(\Z,0,0)$ (with $0$ as
characteristic form). We have $G_{0} = \Q/\Z$, $\phi_{0,0}([x]) =
0$. Trivially any automorphism of $\Q/\Z$ is an automorphism of
the degenerate quadratic function $(\Q/\Z,0)$. However, not every
automorphism of $\Q/\Z$ lifts to an automorphism of $\Z$.
For each bilinear lattice $(M,f)$, Lemma \ref{lem:kernels} gives an
isomorphism $\langle -, - \rangle_{\centerdot}: G_{f} \to
\hbox{Hom}\left(\hbox{Coker}\ \widehat{f}, \Q/\Z\right)$. Therefore
any isomorphism $\psi: \hbox{Coker}\ \widehat{f} \to \hbox{Coker}\
\widehat{f'}$ induces an isomorphism $\psi^{\sharp}:G_{f'} \to
G_{f}$. The map
$$\psi \mapsto \psi^{\sharp}, \ \ \hbox{Iso}\left(\hbox{Coker}\ \widehat{f},
\hbox{Coker}\ \widehat{f'}\right)
\to  \hbox{Iso}\left(G_{f'}, G_{f}\right)$$
is injective (but not surjective in general :
see Lemma \ref{lem:equivalence-utile} below).
We now state the main theorem of this section.

\begin{theorem} \label{th:stabilisation}

Two bilinear lattices $(M,f,c)$ and $(M',f',c')$ with characteristic
forms are stably equivalent if and only if there exists  an element
$$\psi^{\sharp} \in
{\rm{Im}} \left({\rm{Iso}}\left({\rm{Coker}}\ \widehat{f},
{\rm{Coker}}\ \widehat{f}'\right) \to
{\rm{Iso}}\left(G_{f'}, G_{f}\right) \right)$$
such that the associated
quadratic functions $(G_f,\phi_{f,c})$ and $(G_{f'},\phi_{f',c'})$
are isomorphic via $\psi^{\sharp}$.
Furthermore, any such isomorphism  between $(G_{f'},\phi_{f',c'})$ and
$(G_{f},\phi_{f,c})$ lifts to a stable equivalence
between $(M,f,c)$ and $(M',f',c')$.
\end{theorem}
In the sequel, $(M,f,c)$ and $(M',f',c')$ denote bilinear lattices with
characteristic forms. Set $N_{f} = \hbox{Ker}\ \widehat{L}_f$
and  $N_{f'} = \hbox{Ker}\ \widehat{L}_{f'}$.
\begin{lemma} \label{lem:degenerate-part}
Any isomorphism $\psi: {\rm{Coker}}\ \widehat{f} \to
{\rm{Coker}}\ \widehat{f'}$ induces an isomorphism $
{\rm{Ker}}\ \widehat{f'} \to {\rm{Ker}}\ \widehat{f}$. Furthermore,
if $\psi([c]) = [c']$ then $\psi$ induces an isomorphism from the triple
$\left({\rm{Ker}}\ \widehat{f'},0,c'|_{{\rm{Ker}}\ \widehat{f}'}\right)$
onto the triple
$\left({\rm{Ker}}\ \widehat{f},0,c|_{{\rm{Ker}}\ \widehat{f}}\right)$.
\end{lemma}
\begin{proof}
Consider the following diagram:
$$\xymatrix{
0 \ar[r] & {\hbox{T}}\ \hbox{Coker}\ \widehat{f} \ar[r]
\ar[d]_-{\psi|_{ {\textrm{T}}{\textrm{Coker}}\ \widehat{f}}}^{\simeq} &
\hbox{Coker}\ \widehat{f} \ar[r] \ar[d]^{\psi}_{\simeq} &
\left(\hbox{Ker}\ \widehat{f}\right)^{*} \ar[r] \ar[d]^{[\psi]} & 0\\
0 \ar[r] & {\hbox{T}}\ \hbox{Coker}\ \widehat{f'} \ar[r] & \hbox{Coker}\
\widehat{f'} \ar[r] &  \left(\hbox{Ker}\ \widehat{f'}\right)^{*} \ar[r] & 0\\
}$$ where the two exact rows are given by Lemma
\ref{lem:cokernel=ker*} and where $[\psi]$ is induced by $\psi$
and is an isomorphism. Applying $\hbox{Hom}(-,\Z)$ to $[\psi]$
yields an isomorphism $\left(\hbox{Ker}\ \widehat{f'}\right)^{**}
= \hbox{Ker}\ \widehat{f'} \to \left(\hbox{Ker}\
\widehat{f}\right)^{**} = \hbox{Ker}\ \widehat{f}$. For the second
statement,  the image of $[c] \in \textrm{Coker}\ \widehat{f}$ in
$(\hbox{Ker}\ \widehat{f})^{*}$ is just $c|_{{\textrm{Ker}}\
\widehat{f}}$.
 Thus $\psi([c]) = [c']$ implies that $[\psi](c|_{{\rm{Ker}}\ \widehat{f}}) =
c'|_{{\rm{Ker}}\ \widehat{f}'}$, as desired.
\end{proof}
\begin{lemma} \label{lem:equivalence-utile}
Let $\Psi \in {\rm{Iso}}(G_{f'},G_{f})$. Then,
the following assertions are equivalent:
\begin{enumerate}
\item[$(1)$] there exists $\psi \in  {\rm{Iso}}\left({\rm{Coker}}\ \widehat{f},
{\rm{Coker}}\ \widehat{f'}\right)$ such that $\Psi = \psi^{\sharp}$;
\item[$(2)$] $\Psi\left(N_{f'}\right)=N_f$ and the
map $\Psi|_{N_{f'}}:N_{f'} \to N_{f}$ lifts to an isomorphism
${\rm{Ker}}\ \widehat{f'} \to {\rm{Ker}}\ \widehat{f}$.
\end{enumerate}
\end{lemma}
\begin{proof}
$(1) \Longrightarrow (2)$: Lemma \ref{lem:degenerate-part} gives
an isomorphism  ${\rm{Hom}}([\psi],\Z): {\rm{Ker}}\ \widehat{f'}
\to {\rm{Ker}}\ \widehat{f}$. Since $\Psi = \psi^{\sharp}$ corresponds to
 ${\rm{Hom}}(\psi,\Q/\Z)$ via the right adjoint map
$\langle - , - \rangle_{\centerdot}$, ${\rm{Hom}}([\psi],\Z)$ is a
lift of $\Psi|_{N_{f'}}:N_{f'} \to
N_{f}$. \\

$(2) \Longrightarrow (1)$: for $[\alpha] \in {\rm{Coker}}\
\widehat{f}$, consider the homomorphism $$h: [x] \mapsto \langle
[\alpha], \Psi([x]) \rangle, \ \ G_{f'} \to \Q/\Z.$$ Since by
hypothesis, $\Psi|_{N_{f'}}$ lifts to Ker $\widehat{f}'$, the map
$h|_{N_{f'}}$ also lifts to a homomorphism Ker $\widehat{f'} \to \Z$.
Hence by Corollary \ref{co:which-homomorphisms} applied to
$h|_{N_{f'}}$, we obtain that $h:G_{f'} \to \Q/\Z$ lies in the
image of $\hphantom{a}_{\centerdot}\langle -, - \rangle :
{\hbox{Coker}}\ \widehat{f}' \to {\rm{Hom}}(G_{f'},\Q/\Z)$. Thus
there exists $[\beta] \in {\rm{Coker}}\ \widehat{f}'$ such that $h
= \langle [\beta], - \rangle$. In other words, $$\forall [x] \in
G_{f'}, \ \ \  \langle [\alpha], \Psi([x]) \rangle = \langle
 [\beta], [x] \rangle  \in \Q/\Z.$$ The assignment $\psi: [\alpha] \mapsto
[\beta]$ is
 additive and is bijective since $\Psi$ is bijective. By construction, $\Psi =
 \psi^{\sharp}$.
\end{proof}
\begin{proof}[Proof of Theorem \ref{th:stabilisation}]
The nondegenerate case is treated by
Proposition \ref{prop:stabilisation-nondegeneree}.
Consider now the general case.  Let ${\mathfrak{s}}$ be a stable equivalence between
lattices, that is, an isomorphism of symmetric bilinear lattices
${\mathfrak{s}}:(M',f',c') \oplus (U',g',u') \to (M,f,c) \oplus (U,g,u)$
where $(U',g',u')$ and $(U,g,u)$ are unimodular lattices equipped
with characteristic forms.
The map ${\mathfrak{s}}$ induces via the discriminant construction an isomorphism
$\Psi: G_{f'} \to G_f$ between $\phi_{f',c'}$ and $\phi_{f,c}$
since unimodular lattices induce trivial quadratic functions.
The map $s$ also induces in the obvious way an isomorphism
$\psi:\hbox{Coker}\ \widehat{f} \to \hbox{Coker}\ \widehat{f'}$, and it
is easily verified that $\psi^\sharp=\Psi$.\\

Conversely, suppose that the given isomorphism $\Psi$ between $(G_{f'},
\phi_{f',c'})$
and $(G_{f},\phi_{f,c})$ is induced by an isomorphism
$\psi:{\rm{Coker}}\ \widehat{f} \to {\rm{Coker}}\ \widehat{f}'$.
First, the homogeneity defects are preserved by $\Psi$:
$d_{\phi_{f',c'}} = d_{\phi_{f,c}} \circ \Psi$. Since $\Psi = \psi^{\sharp}$,
it follows from (\ref{eq:homogeneity-defect}) that $\psi([c]) = [c']$.
By Lemma \ref{lem:degenerate-part}, the induced map
${\rm{Ker}}\ \widehat{f'}\to
{\rm{Ker}}\ \widehat{f}$ is an isomorphism between
$\left({\rm{Ker}}\ \widehat{f'},0,c'|_{{\rm{Ker}}\ \widehat{f}'}\right)$ and
$\left({\rm{Ker}}\ \widehat{f},0,c|_{{\rm{Ker}}\ \widehat{f}}\right)$.

Let $\left(\bar{M},\bar{f}\right)$ be the nondegenerate lattice
induced by $f$ (see (\ref{eq:splitting_of_M})). Choose  a section
$s$ of the canonical projection $p:M \to \bar{M}$. Lemma
\ref{lem:retract} yields a characteristic form $\bar{c}=
s^{*}|_{{\textrm{Char}}(f)}(c) \in \hbox{Char}(\bar{f})$. The
section $s$ of $p$ induces a map $s:G_{\bar{f}}\to G_{f}$ such
that $\phi_{\bar{f},\bar{c}} = \phi_{f,c} \circ s$.
 Then $s' = \Psi^{-1} \circ
s \circ (\psi|_{{\rm{T}}({\rm{Coker}}\ \widehat{f})})^{-1}$ is a
section of $B_{f'}:G_{f'} \to  {\rm{T}}({\rm{Coker}}\
\widehat{f}')$. Note that by Lemma \ref{lem:sections}, $s'$
is induced by a section of $p':M' \to \bar{M'}$, again denoted $s'$. If
we set $\bar{c'}= (s')^{*}|_{{\textrm{Char}}(f')}(c') \in
\hbox{Char}(\bar{f'})$, we find that $\phi_{\bar{f'},\bar{c'}} =
\phi_{f',c'} \circ s'$. It follows from the definitions that
$\phi_{\bar{f'},{\bar{c'}}} \circ \psi|_{{\rm{T}}({\rm{Coker}}\
\widehat{f})} = \phi_{\bar{f},{\bar{c}}}$. Thus
$\phi_{\bar{f'},{\bar{c'}}} \simeq \phi_{\bar{f},{\bar{c}}}$.
Since these quadratic functions are nondegenerate,  Proposition
\ref{prop:stabilisation-nondegeneree} applies: the isomorphism
$\psi|_{{\rm{T}}({\rm{Coker}}\ \widehat{f})}$ lifts to a stable
equivalence between the lattices $(\bar{M},\bar{f},\bar{c})$ and
$(\bar{M}', \bar{f}',\bar{c}')$.

Finally, there is a stable equivalence between $(M,f,c) \simeq (\bar{M},\bar{f},\bar{c})
 \oplus (\hbox{Ker}\ \widehat{f},0,c|_{{\rm{Ker}}\ \widehat{f}})$ and
 $(M',f',c') \simeq (\bar{M}',\bar{f}',\bar{c}') \oplus (\hbox{Ker}\
 \widehat{f}',0,c'|_{{\rm{Ker}}\ \widehat{f}'})$ where the isomorphisms are induced
by the sections $s$ and $s'$ respectively.
By construction, the isomorphism $G_{f'}\to G_f$ induced by this stable
equivalence is  $\Psi$.
\end{proof}
Consider two bilinear pairings, denoted $\pm 1$, defined on $\Z$ by
$(1,1) \mapsto \pm 1$, and both equipped with the Wu class $1 \in
\Z$. For $n \in \N$, we denote by $n(\Z,\pm 1,1)$ an $n$-fold
orthogonal sum of some copies of $(\Z,\pm 1,1)$. The next result
says that the stable equivalence in Theorem \ref{th:stabilisation}
can be realized by adding orthogonal summands of $(\Z, 1, 1)$ and
$(\Z, -1, 1)$.

\begin{corollary}\label{cor:stabsimple}
Let $(M,f,c)$ and $(M',f',c')$ be two symmetric bilinear lattices equipped
with characteristic forms. We have: $(G_{f},\phi_{f,c})
\simeq (G_{f'},\phi_{f',c'})$ if and only if
there exists $n, n' \in \N$ such that $(M,f,c) \oplus n(\Z,\pm 1,1)
\simeq (M',f',c') \oplus n'(\Z,\pm 1,1)$.
\end{corollary}
%
%
%
%
\begin{proof}
One direction is obvious. For the converse, we apply
Theorem \ref{th:stabilisation}. We obtain unimodular lattices $(U,g)$
and $(U',g')$, equipped with characteristic forms $s$ and $s'$ respectively,
such that there is an isomorphism  sending $(M,f,c) \oplus (U,g,u)$ onto
$(M',f',c') \oplus (U',g',u')$. By stabilizing if necessary with $(\Z,1)$ and
$(\Z,-1)$, we may assume that, as a symmetric bilinear pairing, $(U,g)$
is indefinite\footnote{i.e. $g$ takes both positive and negative values.}
and not even\footnote{even means that $g(x,x) \in 2\Zi$ for all $x \in U$.}.
Then a theorem asserts that $(U,g)$ is isomorphic to an orthogonal
sum of copies of $(\Z, 1)$ and $(\Z, -1)$ \cite{Ser}. It follows that $(U,g,u)$
is isomorphic to a sum of copies of $(\Z, 1, 1)$ and $(\Z, -1, 1)$.
(Here we use the fact that any odd integer $a \in \Z$ induces a characteristic
form for $(\Z, \pm 1)$ and any such triple $(\Z,\pm 1,a)$ is isomorphic to
$(\Z, \pm 1, 1)$. More generally, it follows from
\S \ref{subsec:latticesandforms} that any unimodular bilinear lattice
$(U,g)$ has only one characteristic form modulo $2\widehat{g} (U)$.) A similar observation holds for $(U',g',u')$.
We conclude that there is a stable equivalence between
$(M,f,c)$ and $(M',f',c')$ involving only a stabilization with copies of
$(\Z, 1, 1)$ and $(\Z, -1, 1)$, which is the desired result.
\end{proof}
\begin{remark}[and erratum]
This is a generalization of \cite[Lemma 2.1(b)]{De1} (whose proof
is unfortunately incorrect).
\end{remark}
\section{A complete system of invariants}
\label{sec:system}

We present a complete system of invariants for a certain class of
quadratic functions on torsion abelian groups, including the
nondegenerate quadratic functions defined on finite abelian groups.

First, we assume that the abelian group $G$ is finite. Let ${S}^{1}$
be the multiplicative group of complex numbers of absolute value
$1$.  For any quadratic function $q:G \to \Q/\Z$, we define the
Gauss sum
$$\gamma(q) = |G|^{-\frac{1}{2}} \sum_{x \in G} \exp(2\pi iq(x)) \in \C.$$
It is not difficult to see that $\gamma(q) \in S^1$ if $q$ is
nondegenerate.  An important observation is  the relation
\begin{equation}
\gamma(\alpha \cdot q) = e^{-2\pi iq(\alpha)} \gamma(q), \ \alpha
\in G. \label{eq:Gaussrel}
\end{equation}
Recall that $b_q$ (resp. $d_q=q-\overline{q}\in G^*$) is the
bilinear pairing (resp.  the homogeneity defect) associated to
$q$. The following result is the main result of this section.
\begin{theorem}
\label{th:invcomplet}  Two nondegenerate quadratic functions $q:G
\to \Q/\Z$ and $q':G' \to \Q/\Z$ on finite abelian groups are
isomorphic if and only if there is an isomorphism $\psi:G \to G'$
such that $\psi^{*}b_{q'} = b_{q}$, $\psi^{*}d_{q'} = d_{q}$ and
$\gamma(q') = \gamma(q)$.
\end{theorem}

\noindent That two isomorphic quadratic functions have same Gauss
sums and isomorphic associated bilinear
pairings and homogeneity defects is straightforward. The
nonobvious part lies in the converse.

Since homogeneous quadratic functions have trivial homogeneity
defect, we deduce the following result (see \cite[Theorem
1.11.3]{Nik}).

\begin{corollary}[V. Nikulin] \label{th:homogeneous}
Two nondegenerate homogeneous quadratic functions $q:G\to \Q/\Z$ and
$q':G'\to \Q/\Z$ on finite abelian groups are isomorphic if and only
if there is an isomorphism $\psi:G\to G'$ such that $\psi^{*}b_{q'}
= b_{q}$ and $\gamma(q') = \gamma(q)$.
\end{corollary}

The key step in the proof of Theorem \ref{th:invcomplet} is the
following fundamental lemma.
\begin{lemma}[Fundamental Lemma]\label{lem:fundamental}
Let $q:G\to \Q/\Z$ be a quadratic function on a torsion abelian
group $G$. If $\alpha \in G$ is an element of order $2$ such that
$q(\alpha) = 0$, then $\alpha \cdot q \sim q$.
\end{lemma}
\begin{proof}[Proof of Th. \ref{th:invcomplet}
from Lemma \ref{lem:fundamental}] With no loss of generality, we
may assume that $G = G'$, $b_{q} = b_{q'}$, $d_{q} = d_{q'}$ and
$\gamma(q)=\gamma(q')$ and then show that $q \sim q'$. The
equality $b_{q} = b_{q'}$ implies (cf. \S \ref{sec:quad}) that $q'
= \alpha \cdot q$ for some $\alpha \in G$. The equality $d_{q} =
d_{q'}$ implies $2\alpha = 0$. The equality $\gamma(q) =
\gamma(q')$, together with $(\ref{eq:Gaussrel})$ imply that
$q(\alpha) = 0$. Therefore Lemma \ref{lem:fundamental} applies.
\end{proof}

Hence it remains to prove Lemma \ref{lem:fundamental}.

\begin{proof}[Proof of the Fundamental Lemma]
If we were just looking for a permutation $\sigma$ of $G$ such
that $q\circ \sigma=\alpha \cdot q$, we could simply choose the
involution $\left(x \mapsto x + \alpha\right)$ since $q(x+\alpha)
= q(x) + b_{q}(x,\alpha) + q(\alpha) = q(x) +
b_{q}(\alpha,x)=(\alpha \cdot q )(x)$, $x \in G$. So it is
sufficient to find $\phi\in \hbox{End}\left(G\right)$ such that
\begin{equation}
(\hbox{Id}_{G} + \phi)^{2} = {\rm{Id}}_{G} \ \ {\rm{and}}\ \  q(x +
\phi(x)) = q(x + \alpha), \ {\hbox{for any}}\ x\in G.
\label{eq:suffices}
\end{equation}
Since $\alpha$ has order $2$, the map $\widehat{b}_{q}(\alpha) \in
G^{*}$ is of order 2. For any $x \in G$, define $n(x) \in \Z/2\Z$ by
$$\frac{n(x)}{2} = b_{q}(\alpha, x) \in \Q/\Z.$$
Consider the map $\phi:G \to G$ defined by
$$\phi(x) = n(x)\alpha, \ \ x \in G.$$ It is a well-defined
(independent of the lift in $\Z$ of $n(x) \in \Z/2\Z$)
endomorphism of $G$. Using the fact that $q(\alpha) = 0$, we
 readily verify that $\phi$ satisfies
(\ref{eq:suffices}).
\end{proof}

We now deduce a more general version of Theorem \ref{th:invcomplet}
which allows for some degeneracy.
 Consider the following commutative diagram of extensions of abelian
 groups
 $$\xymatrix{
 0 \ar[r] & A \ar[r]^{i} \ar[d]_{\simeq}^{\psi|_{A}} &
 G \ar[d]^{\psi}_{\simeq} \ar[r]^{p} & B \ar[r] \ar[d]_{\simeq}^{[\psi]} & 0 \\
 0 \ar[r] & A' \ar[r]^{i'}  & G' \ar[r]^{p'} & B' \ar[r] & 0.
 }$$
 Here $[\psi]:B \to B'$ is the isomorphism induced by $\psi$.

\begin{definition} \label{def:compatible}
Given a diagram as above, two sections $s$ and $s'$ of $p$ and $p'$
respectively are said $\psi$-{\emph{compatible}}, or simply
{\emph{compatible}}, if $s' = \psi \circ s \circ [\psi]^{-1}$.
$$\xymatrix{
G \ar[r]_{p} \ar[d]_{\psi} & B \ar[d]^{[\psi]} \ar@/_/@<-1ex>[l]_{s} \\
G' \ar[r]^{p'} & B' \ar@/^/@<1ex>[l]^{s'}
}$$
\end{definition}
Let $q: G \to \Q/\Z$ be a quadratic function on a torsion abelian
group. Clearly $G$ can be regarded as the extension of $A =$ Ker
$\widehat{b}_q$ by $B = G/\hbox{Ker}\ \widehat{b}_q$. We shall say
that $(G,q)$ meets the {\emph{finiteness condition}} if the
following two conditions are satisfied:
\begin{enumerate}
\item[-] $G/\hbox{Ker}\ \widehat{b}_{q}$ is finite;
\item[-] the extension $G$ of $\hbox{Ker}\ \widehat{b}_{q}$ by
$G/\hbox{Ker}\ \widehat{b}_{q}$ is split.
\end{enumerate}
\begin{example}
Let $(M,f,c)$ be a symmetric bilinear lattice
equipped with a characteristic form $c$.
Then, the quadratic function $\left(G_f,\phi_{f,c}\right)$
meets the finiteness condition since by (\ref{eq:splitq}), the short exact sequence (\ref{eq:shortexact}) is split.
\end{example}

For any given section $s$
of the projection $p:G \to G/{\rm{Ker}}\ \widehat{b}_{q}$, $q \circ s$ is a
nondegenerate quadratic
function on $G/{\rm{Ker}}\ \widehat{b}_{q}$. Note that the finiteness
condition implies that $\gamma(q \circ s)$ is well-defined for any section $s$.
In the sequel, we denote by $r_q: \hbox{Ker}\ \widehat{b}_q\to \Q/\Z$ the
homomorphism $q|_{\hbox{\scriptsize{Ker}}\ \widehat{b}_q}$.
\begin{corollary} \label{cor:invcomplet-deg}
Two quadratic functions $q:G\to \Q/\Z$ and $q':G'\to \Q/\Z$ satisfying the finiteness condition  are
isomorphic if and only if there is an isomorphism $\psi:G\to G'$ such that
$\psi^{*}b_{q'} = b_{q}$, $\psi^{*}d_{q'} = d_{q}$, $\psi^{*}r_{q'} = r_{q}$
and $\gamma(q' \circ s') = \gamma(q \circ s)$
for $\psi$-compatible sections $s$ and $s'$.
\end{corollary}
\begin{proof} Consider the nondegenerate quadratic functions $q_1=q\circ s$ and $q_1'=q'\circ s'$. The isomorphism $\psi$ induces an isomorphism $[\psi]: G/\hbox{Ker}\ \widehat{b}_q \to G'/\hbox{Ker}\ \widehat{b}_{q'}$. It is readily checked that $[\psi]^{*}b_{q_1'} = b_{q_1}$,
$[\psi]^{*}d_{q_1'} = d_{q_1}$ and $\gamma(q_1)=\gamma(q'_1)$.
By Theorem \ref{th:invcomplet},
$q_1$ and $q'_1$ are isomorphic. We deduce that
$\left(G,q\right)\simeq \left(G/\hbox{Ker}\ \widehat{b}_q,q_1\right)
\oplus \left(\hbox{Ker}\ \widehat{b}_{q},r_{q}\right)$ is isomorphic to
$\left(G',q'\right)\simeq \left(G'/\hbox{Ker}\ \widehat{b}_{q'},q'_1\right)
\oplus \left(\hbox{Ker}\ \widehat{b}_{q'},r_{q'}\right)$ .
\end{proof}

\begin{remark} \label{rem:nb}
Theorem \ref{th:invcomplet} does not say that if $\psi^{*}b_{q'} =
b_{q}$, $\psi^{*}d_{q'} = d_{q}$ and $\gamma(q') = \gamma(q)$,
then $\psi^{*}q' = q$. In general, the isomorphism between $q$ and
$q'$ will be different from $\psi$. However, how the isomorphism
between $q$ and $q'$ will differ from $\psi$ can be read off from
the proof of the Fundamental Lemma.
\end{remark}

\bibliographystyle{amsalpha}

\end{document}